\documentclass[numbers,webpdf]{ima-authoring-template}%

\graphicspath{{Fig/}}

\theoremstyle{thmstyletwo}%
\newtheorem{theorem}{Theorem}
\newtheorem{proposition}[theorem]{Proposition}%

\newtheorem{definition}{Definition}

\numberwithin{equation}{section}
\usepackage{algorithm}
\usepackage{algpseudocode}
\usepackage{amsmath}
\usepackage{amssymb}
\usepackage{array}
\usepackage{breqn}
\usepackage{graphicx}
\usepackage{hyperref}
\usepackage[utf8]{inputenc}
\usepackage{xfrac}
\usepackage{xcolor}
\definecolor{mygreen}{rgb}{0,0.6,0}
\definecolor{mygray}{rgb}{0.5,0.5,0.5}
\definecolor{mymauve}{rgb}{0.58,0,0.82}
\definecolor{altblue}{rgb}{0.0,0.6,1.0}
\definecolor{lstbg}{cmyk}{0.05, 0.01, 0, 0}
\definecolor{morebluish}{cmyk}{0.06,0.04,0,0}
\usepackage{listings}
\input{listings-maple-definition.sty}
\usepackage{mathtools}
\DeclarePairedDelimiter{\ceil}{\lceil}{\rceil}
\DeclarePairedDelimiter{\floor}{\lfloor}{\rfloor}
\usepackage{subfigure}

\newcommand{\A}{\ensuremath{\mathbf{A}}}
\newcommand{\B}{\ensuremath{\mathbf{B}}}
\newcommand{\C}{\ensuremath{\mathbf{C}}}
\newcommand{\D}{\ensuremath{\mathbf{D}}}
\newcommand{\E}{\ensuremath{\mathbf{E}}}
\newcommand{\G}{\ensuremath{\mathbf{G}}}
\newcommand{\I}{\ensuremath{\mathbf{I}}}
\newcommand{\J}{\ensuremath{\mathbf{J}}}
\renewcommand{\P}{\ensuremath{\mathbf{P}}}

\renewcommand{\S}{\ensuremath{\mathbf{S}}}

\renewcommand{\emph}[1]{\textsl{#1}}
\renewcommand{\textit}[1]{\emph{#1}}

\newcommand{\Char}{\ensuremath{\mathcal{C}}}
\newcommand{\Fam}{\ensuremath{\mathcal{B}}}

\usepackage{xfrac}
\usepackage{tabularx,booktabs}

\newcommand{\res}{\ensuremath{R}}
\newcommand{\Sc}{\ensuremath{\mathbf{S}_C}}  

\def\sp{\mathrm{spread}}
\def\SP{\mathrm{MaxSpread}}

\begin{document}

\DOI{DOI HERE}
\copyrightyear{2025}
\vol{00}
\pubyear{2025}
\access{Advance Access Publication Date: Day Month Year}
\appnotes{Paper}
\copyrightstatement{}
\firstpage{1}


\title[Symmetric Bohemian Spread]{On the maximal spread of symmetric Bohemian matrices}

\author{Neil J.~Calkin\ORCID{0000-0002-5701-8581}
\address{\orgdiv{Dept. of Mathematics}, \orgname{Clemson University}, \orgaddress{\street{105 Sikes Hall}, \postcode{29634}, \state{South Carolina}, \country{USA}}}}
\author{Robert M.~Corless\ORCID{0000-0003-0515-1572}
\address{\orgdiv{Dept. of Computer Science}, \orgname{Western University}, \orgaddress{\street{1151 Richmond Street}, \postcode{N6A 3K7}, \state{Ontario}, \country{Canada}}}}
\author{Laureano Gonzalez-Vega\ORCID{0000-0002-3934-3890}
\address{\orgdiv{Dept. de Métodos Cuantitavos}, \orgname{CUNEF Universidad}, \orgaddress{\street{C. de Almansa, 101}, \postcode{28040}, \state{Madrid}, \country{Spain}}}}
\author{J. Rafael Sendra\ORCID{0000-0003-2568-1159}
\address{\orgdiv{Dept. de Matemáticas}, \orgname{CUNEF Universidad}, \orgaddress{\street{C. de Almansa, 101}, \postcode{28040}, \state{Madrid}, \country{Spain}}}}
\author{Juana Sendra\ORCID{0000-0002-9927-169X}
\address{\orgdiv{Dept. de Matemáticas}, \orgname{CUNEF Universidad}, \orgaddress{\street{C. de Almansa, 101}, \postcode{28040}, \state{Madrid}, \country{Spain}}}}

\authormark{Calkin et al.}

\corresp[*]{Corresponding author: Robert M.~Corless \href{email:rcorless@uwo.ca}{rcorless@uwo.ca}}

\received{Date}{0}{Year}
\revised{Date}{0}{Year}
\accepted{Date}{0}{Year}


\abstract{Let $\A$ be a square matrix with entries in $\mathbb{R}$. The spread of $\A$ is defined as the maximum of the distances among the eigenvalues of $\A$. Let $S_m[a,b]$ denote the set of all $m\times m$ symmetric matrices with entries in the real interval $[a,b]$ and let $S_m\{a,b\}$ be the subset of $S_m[a,b]$ of Bohemian matrices with population from only the extremal elements $\{a,b\}$. S.~M.~Fallat and J.~J.~Xing in 2012 (see~\cite{fallat2012spread}) proposed the following conjecture:  the maximum spread in $S_m[a,b]$ is attained by a rank $2$ matrix in $S_m\{a,b\}$. X.~Zhan had proved previously that the conjecture was true for $S_m[-a,a]$ with $a>0$. We will show how to interpret this problem geometrically, via polynomial resultants, in order to be able to treat this conjecture from a computational point of view. This will allow us to prove that this conjecture is true for several formerly open cases.}


\maketitle

\begin{center}
    \textbf{Dedicated to the memory of Nicholas J. Higham}
\end{center}
\section{Introduction}


Let $\A$ be an $m\times m$ square matrix with entries in $\mathbb{C}$. The spread of $\A$ is defined to be
\[ \sp(\A)=\max\{|\lambda_i-\lambda_j|  \colon i,j\in \{1,\ldots,m\}\}, \]
where $\lambda_1,\ldots,\lambda_m$ are the eigenvalues of $\A$.
In this paper, we   study  the spread of a matrix in the class of the symmetric $m\times m$ matrices with entries in a real interval $[a,b]$. We denote this set as  $S_m([a,b])$ and let $S_m(\{a,b\})$ denote the set of symmetric $m\times m$ matrices with entries in the finite discrete set $\{a,b\}$.

The question is to determine 
\[ \SP(S_m([a,b])):=\max\{ \sp(\A)\colon \A\in S_m([a,b])\}. \]
X.~Zhan asked in~\cite{Zhan2005} (see Problem 2) to compute $\SP(S_m([a,b]))$ and to determine which matrices attain the maximum. He solved this problem for symmetric intervals: if $\A\in S_m([-a,a])$ with $m\geq 2$ then
\begin{equation} 
\sp(\A)\leq
\left\{\begin{array}{ll}
a\sqrt{2m^2}&\hbox{if $m$ is even}\\
\noalign{\vspace*{1mm}}
a\sqrt{2m^2-1}& \hbox{if $m$ is odd}
\end{array}\right. \label{eq:Zhanmax1}
\end{equation}
Let $\J_{r,s}$ denote  the $r \times s$ matrix with all entries equal to $1$, and  $\J_r=\J_{r,r}$. Then,   
if $m$ is even, the equality in \eqref{eq:Zhanmax1}  holds if and only if $\A$ is sign permutation similar to
$$a\left(\begin{array}{cc}
1&1\\
1&-1\end{array}\right)\otimes \J_{\frac{m}{2}} \>.$$
If $m$ is odd, the equality in \eqref{eq:Zhanmax1} holds iff $\A$ is sign permutation similar to
$$a\left(\begin{array}{cc}
\J_{\frac{m+1}{2}}&\J_{\frac{m+1}{2},\frac{m-1}{2}}\\
\noalign{\vspace*{1mm}}
\J_{\frac{m-1}{2},\frac{m+1}{2}}&-\J_{\frac{m-1}{2}}\end{array}\right)\>.$$

Later Fallat and Xing introduced in~\cite{fallat2012spread} the matrix candidate to provide the sought maximum spread in $S_m([a,b])$:
\begin{quotation} 
\phantomsection\label{conjecture}
  \textsl{The 2012 Fallat and Xing  Conjecture}: for $\A\in S_m([a,b])$, the maximum spread in this class is attained by some $\A\in S_m(\{a,b\})$ with $\mathrm{rank}(\A)=2$.
\end{quotation}

They showed that it is enough to study the case $S_m([a,1])$ with $-1\leq a <1$. They proved moreover that if $\A\in S_m([a,1])$, $\mathrm{rank}(\A)=2$, and $\sp(\A)$ is maximal in this set of matrices, then $\A$ is of the form (up to permutation similarity)
$$\A=
\left(\begin{array}{cc}
a\,\J_k&\J_{k,m-k}\\
\noalign{\vspace*{1mm}}
\J_{m-k,k}&\J_{m-k,m-k}
\end{array}\right)$$
where the blocksize $k$ is given specifically by
$$k=\mathrm{round}\left[\dfrac{m}{a+3}\right]$$
and 
$$
\sp(\A)=\sqrt{(a^2+2a-3) k^2+2m (1-a) k+m^2}\>.$$

Finally Biborski simplified the conjecture in~\cite{biborski2022note} by showing that the maximal spread would occur only when the matrix had entries belonging to the set $\{a,b\}$: $\SP(S_m([a,b])$ can only be attained at $S_m(\{a,b\})$.  That is, only the extremal points of the interval $[a,b]$ were needed. Moreover, by using this result, he solved the conjecture for the case $S_m([0,1])$ with $m\leq 3$ and opened the way to investigate what happens with small values of $m$ for $[0,1]$ and $[a,1]$. \textcolor{black}{This reduction of the problem, from $S_m([a,b])$ to $S_m(\{a,b\})$, shows that the  analysis of the conjecture can be restricted to the set of $m\times m$ Bohemian matrices with population $\{a,b\}$\footnote{A set of matrices with entries from a fixed finite population, usually a subset of $\mathbb{Z}$, is called “Bohemian”. The mnemonic comes from BOunded HEight Matrix of Integers, BOHEMI; see \cite{Thornton2019Bohemian,LMS2020,ChanEtAl2020UHTB} for further details.}.}

In this paper we analyze the \textsl{2012 Fallat and Xing  Conjecture} from a computational point of view. In this way we introduce several algorithms to determine  $\sp(\A)$ and to compute exactly $\SP(S_m([a,1]))$. We have been able to prove the conjecture for the following cases  for $S_m([a,1])$:
\begin{enumerate}
\item When $a=0$ and $m \in \{2, 3, 4, 5, 6, 7, 8\}$ (see Section \ref{sec-casezero}).
\item When $-1 < a < 1$ and $m \in \{2, 3, 4, 5, 6, 7\}$ (see Section \ref{sec:symbolica}).
\item \textcolor{black}{When $a=0$ and $m$ is a multiple of $3$ (see Subsection \ref{subsec-multiple}).}
\end{enumerate}
In order to be able to deal with the huge number of matrices to consider (the number of matrices in $S_m(\{a,1\})$ is $2^{m(m+1)/2}$), several economies have been introduced involving graph techniques allowing us to discard many matrices known not to make the spread bigger than the one previously computed.

Moreover, in all these cases, we have verified that no matrix with $\mathrm{rank}(\A)\neq 2$ attains the maximum spread in $S_m(\{a,1\})$. A computational approach was also followed by Xing in his PhD Thesis (see the Appendix in \cite{Xing}) in order to validate (or not) this conjecture by generating randomly many matrices but he did not pursue the exhaustive analysis we have used here.


As a byproduct of this approach we have been able to prove the conjecture for $S_m([0,1])$ when $m$ is a multiple of $3$.

\textcolor{black}{The paper is structured as follows. In Section \ref{section-numertical}, we discuss various aspects of the spread, supported by numerical evidence indicating that the spreads seems to be approximately normally distributed. Section \ref{section-determining}  investigates different methods for computing the spread and revisit some previous results in this topic. Leveraging tools from graph theory, we prove the conjecture for the case m divisible by $3$ and within the interval $ [0,1].$  Additionally, we establish a link between the problem and the use of resultants. Much of our argument relies on computational methods using mathematical software; we use Maple. In Section \ref{sec-reliable},  we discuss the feasibility and reliability of such proofs. Section \ref{sec-brute-force} considers the case $m = 4$ with $a = 0$ and symbolic $a$. This particular analysis illuminates the general method presented later. A key computational bottleneck of the spread problem is the exponential growth in the number of matrices as $m$ increases. Since many share the same characteristic polynomial, in Section \ref{sec-graphs} we use graph theory to group them. While not a strict partition, this decomposition significantly reduces the computation. Section  \ref{sec-casezero}  is devoted to prove the conjecture for $2\leq m\leq 8$ for the particular case of $a=0$ and Section \ref{sec:symbolica} establishes the validity of the conjecture for $m\in \{2,\ldots,7\}$ and $a$ in the interval $(-1,1)$. \textcolor{black}{The final section, Section \ref{sec:another}, presents a different approach to the study of the spread, which can be viewed as an alternative approach that may be useful in future, and also as a way of understanding what a \textsl{resultant} is if one has not seen the concept before.} The paper ends with a section with the conclusions.}

\section{Some preliminary numerical computation}\label{section-numertical}
Consider symmetric matrices of dimension $m$ whose entries 
are drawn uniformly at random from the interval $[-1,1]$.  We wrote a Julia program that generated a number of samples of such matrices, and used the \lstinline{LinearAlgebra.eigvals} routine to compute the eigenvalues of each sample matrix.  We believe that this routine detects the attribute ``Symmetric'' and calls the appropriate LAPACK routines~\cite{anderson1999lapack} to generate the real eigenvalues, by first reducing to real symmetric tridiagonal form using DSYTRD and then calling DSTERF. The Julia routine then sorts the eigenvalues in ascending order.  The spread of the sample matrix can therefore be computed simply by \lstinline{ev[m]-ev[1]}. 
We then normalize by the known maximal spread of such matrices, which was proved in~\cite{Zhan2005} to be given by equation~\eqref{eq:Zhanmax1}. 

\begin{figure}
    \centering
    \includegraphics[width=0.5\linewidth]{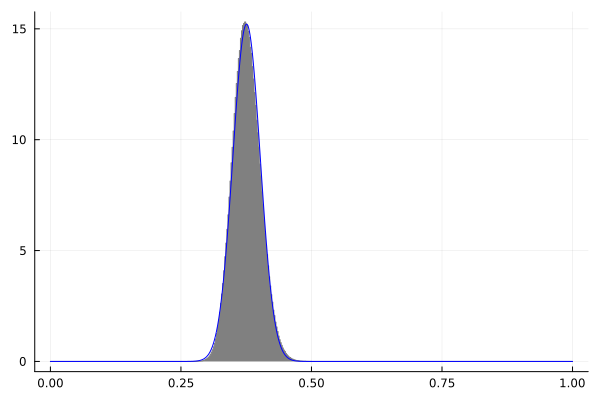}
    \caption{Distribution of spreads of two million symmetric matrices of dimension $m=13$ whose entries were drawn uniformly at random from the interval $[-1,1]$.  Eigenvalues were computed numerically in Julia.  The arithmetic mean of the sampled spreads was $\mu = 0.376$ and the variance was $\sigma^2 = 0.000688$. A normal distribution with those parameters is sketched in blue together with the computed distribution.}
    \label{fig:twomillion}
\end{figure}

We see in figure~\ref{fig:twomillion} that the spreads seem to be distributed normally (we chose $m=13$ and took a sample of two million matrices: the picture did not change appreciably if we took instead twenty million or two hundred million).  This is in apparent agreement with the central limit theorem, in spite of the fact that the maximum normalized spread is just $1$.  

Altering the dimension and running the experiments again gives similar results.  We also see that increasing the dimension $m$ reduces the sample mean as a fraction of the theoretical maximum. Indeed, we note that the arithmetic mean of the unnormalized spreads seems to be $O(1)$ as the dimension $m$ goes to infinity, but we have not done enough experiments for this to be more than a speculation.

This computation is quite informative, and tells us that having a spread close to the theoretical maximum requires quite a rare matrix.  Most sample matrices will have spreads only a small fraction of the maximum possible, unless the dimension is very small.  Nonetheless, the maximum possible spread is important in some situations.

We leave off numerical experiments now and concentrate on the theoretical maximum spread.  It turns out that the maximum is not known for all possible choices of intervals 
$[a,b]$ 
to draw from; symmetric intervals 
$[-a,a]$ as we examined here are the exception, where the maximum spread is known to be $O(a\sqrt{m})$ (a more precise statement occurs in equation~\eqref{eq:Zhanmax}). 

Moreover, numerical computation of eigenvalues will give accurate estimates of the maximum spread only if computation is exhaustive in some fashion, and that this is possible is not immediately evident.  Luckily, as we will see in the next section, we may restrict attention to a subclass of matrices that is finite in number at each dimension.  Somewhat unluckily, however, this finite number grows combinatorially with the dimension, and so computation will only get us so far.

\section{Determining the maximum spread in $S_m([a,b])$}\label{section-determining}
\textcolor{black}{In the introduction, we reviewed some of the previous achievements on the study of the spread. In     this section, we discuss these results in greater detail, we prove the conjecture for $m$ being a multiple of 3 and the interval $[0,1]$, and we provide the idea on how to approach the problem via resultants.}

\subsection{Bounding the spread}
Let $\A$ be an $m\times m$ square matrix with entries in $\mathbb{C}$. We define the spread of $\A$ as
\[ \sp(\A)=\max\{|\lambda_i-\lambda_j|  \colon i,j\in \{1,\ldots,m\}\}, \]
where $\lambda_1,\ldots,\lambda_m$ are the eigenvalues of $\A$.

In 1956 Mirsky stated the following bound (see \cite{Mirsky1956Spread}): for  $m\geq 3$ 
\begin{equation}\label{MirskyBound}
\sp(\A)\leq \sqrt{2\|\A\|_{F}^{2}-\dfrac{2}{m}|\mathrm{Tr}(\A)|^2},
\end{equation}
\textcolor{black}{where $\|\A\|_{F}$ and $\mathrm{Tr}(\A)$ denote the Frobenius norm and the trace of $\A$,  respectively.}
{Moreover}, the equality holds if and only if $\A$ is normal and $m-2$ of its eigenvalues are equal and this common value is the arithmetic mean of the  two others. 

In particular, we focus on the class of the symmetric $m\times m$ matrices with entries in a real interval $[a,b]$. We denote this set as 
$S_m([a,b])$. Let $S_m(\{a,b\})$ denote the set of symmetric $m\times m$ matrices with entries in the finite discrete set $\{a,b\}$.

The question is to analyze/compute 
\[ \SP(S_m([a,b])):=\max\{ \sp(\A)\colon \A\in S_m([a,b])\}. \]

X.~Zhan asked in~\cite{Zhan2005} to compute $\SP(S_m([a,b]))$ but Fallat and Xing introduced in~\cite{fallat2012spread} the matrix candidate to provide the searched maximum spread; \textcolor{black}{see \hyperref[conjecture]{ 2012 Fallat and Xing Conjecture}   in page~\pageref{conjecture}.}


They show in~\cite{fallat2012spread} that it is enough to study the case $S_m([a,1])$ with $-1\leq a <1$.

\begin{theorem} \cite{fallat2012spread}
Let $\A\in S_m([a,1])$, with $-1\leq a <1$ and $\mathrm{rank}(\A)=2$. If $\sp(\A)$ is maximal then  $\A$ is of the form (up to permutation similarity)
$$\A=
\left(\begin{array}{cc}
a\,\J_k&\J_{k,m-k}\\
\J_{m-k,k}&\J_{m-k,m-k}
\end{array}\right)$$
where the blocksize $k$ is given specifically by
$$k=\mathrm{round}\left[\dfrac{m}{a+3}\right]$$
and 
$$
\sp(\A)=\sqrt{(a^2+2a-3) k^2+2m (1-a) k+m^2}\>.$$
\end{theorem}

If $m/(a+3) = \ell + 1/2$ is exactly half-way between integers then each rounding, up or down, will give matrices with the same maximal spread. 

Biborski reduced the conjecture in~\cite{biborski2022note} by showing that the maximal spread would occur when the matrix had entries belonging to the set $\{a,b\}$.  That is, only the extremal points of the interval $[a,b]$ were needed.

\begin{theorem} \cite{biborski2022note}
$\SP(S_m([a,b])$ can only be attained at $S_m(\{a,b\})$.
\end{theorem}

Moreover, Zhan solved the conjecture for symmetric intervals $S_m([-a,a])$ in~\cite{Zhan2005}.

\begin{theorem} \cite{Zhan2005}
Let $\A\in S_m([-a,a])$ with $m\geq 2$. Then
\begin{equation}
\sp(\A)\leq 
\left\{\begin{array}{ll}
a\sqrt{2m^2}&\hbox{if $m$ is even}\\
\noalign{\vspace*{1mm}}
a\sqrt{2m^2-1}&\hbox{if $m$ is odd}
\end{array}\right. \label{eq:Zhanmax}
\end{equation}
If $m$ is even then equality holds if and only if $\A$ is sign permutation similar to
$$a\left(\begin{array}{cc}
1&1\\
1&-1\end{array}\right)\otimes \J_{\frac{m}{2}}\>.$$
If $m$ is odd then equality holds iff $\A$ is sign permutation similar to
$$a\left(\begin{array}{cc}
\J_{\frac{m+1}{2}}&\J_{\frac{m+1}{2},\frac{m-1}{2}}\\
\noalign{\vspace*{1mm}}
\J_{\frac{m-1}{2},\frac{m+1}{2}}&-\J_{\frac{m-1}{2}}\end{array}\right)\>.$$
\end{theorem}

\begin{definition}  \cite{fallat2012spread}
 A real symmetric matrix $\A$ of order $m$ ($m\geq 3$) is a Mirsky matrix if the spread of $\A$ attains the Mirsky bound. 
 \end{definition}
 
 Next we simplify the proof in~\cite{Zhan2005} when $m$ is even and $a=-1$.
 
\begin{theorem} The order $m$ ($m\geq 3$) matrix 
$$\A=
\left(\begin{array}{cc}
a\,\J_k&\J_{k,m-k}\\
\J_{m-k,k}&\J_{m-k,m-k}
\end{array}\right)$$  is a Mirsky matrix if and only if $$m=(1-a)\mathrm{round}\left[\dfrac{m}{a+3}\right]$$ if and only if $a=-1$ and $m$ is even.
\end{theorem}

Let $m$ be even and $a=-1$. Then
$$(1-a)\mathrm{round}\left[\dfrac{m}{a+3}\right]=2\ \mathrm{round}\left[\dfrac{m}{2}\right]=m\>.$$

\subsection{\textcolor{black}{Proof of the conjecture for $S_{3q}([0,1])$}}\label{subsec-multiple}

\textcolor{black}{Now, we deal with the particular case of $S_m([0,1])$ when $m$ is a multiple of $3$. We start with the following result from graph theory.}

\begin{theorem}  \cite{breen2022maximum} Let $\A=(a_{i,j})$ be an $m\times m$ symmetric nonnegative matrix. Then 
$$\sp(\A)\leq \dfrac{2m}{\sqrt{3}} \,\max_{i,j} \{a_{i,j}\}\>.$$
\end{theorem}

Then, for any matrix $\A\in S_m([a,1])$, with $0\leq a <1$, we have  $$\sp(\A)\leq\frac{2m}{\sqrt{3}}.$$ 
Let $$B=
\left(\begin{array}{cc}
a\,\J_k&\J_{k,m-k}\\
\J_{m-k,k}&\J_{m-k,m-k}
\end{array}\right)$$ with $0\leq a <1$.
{\color{black}  If $a=0$ and $m=3q$ then
$$\sp(B)^2={(a^2+2a-3) k^2+2m (1-a) k+m^2}=- 3 k^2+6q k+9q^2$$
with
$$k=\mathrm{round}\left[\dfrac{m}{a+3}\right]=\mathrm{round}\left[\dfrac{3q}{3}\right]=q\>.$$
Therefore
$$\sp(B)=\sqrt{-3 k^2+6q k+9q^2}=\sqrt{12q^2}=2\sqrt{3}q\>.$$

Since for any $\A\in S_{3q}([a,1])$ we have
$$\sp(\A)\leq \frac{2m}{\sqrt{3}}\max_{i,j} a_{i,j}=\frac{6q}{\sqrt{3}}=2\sqrt{3}q$$
we conclude that the maximum spread in $S_m([0,1])$ is attained by the rank $2$ matrix $B$.
Summarizing, we get the following result.

\begin{theorem}\label{theorem-multiple 3}
If $a=0$ and the dimension $m$ is a multiple of $3$ then  \textsl{2012 Fallat and Xing  Conjecture}    is true.
\end{theorem}
}

\subsection{Computing symbolically the spread of a matrix in $S_m([a,b])$.}\label{subsection-resultant}

Let $\A$ be a matrix in $S_m([a,b])$ and $p(\lambda)$ the characteristic polynomial of $\A$ such that $p(\lambda)\not = \lambda^m$ (in this case the spread of $\A$ is $0$). Defining  
$$R(T)={\rm resultant}_{\lambda}(p(\lambda),p(T+\lambda))$$
we have that $R(T)$ is a polynomial in  {$\mathbb{R}[T]$} verifying:
\begin{itemize} 
\item If $\lambda_1, \ldots,\lambda_m$ are the eigenvalues of $\A$ then the roots of $R(T)$ are $\lambda_i-\lambda_j$, $1\leq i,j, \leq m$.
\item $\deg_T(R)=m^2$.
\item \textcolor{black}{The multiplicity of $0$ as root of $R$ is bigger than or equal to $m$ because each $\lambda_i-\lambda_i$ is a root of the resultant. Let $\psi(T)=R(T)/T^m$.  $\deg_T(\psi)=m^2-m$.}
\item \textcolor{black}{$\psi(T)$ is a polynomial in $T^2$ because both $\lambda_i-\lambda_j$ and $\lambda_j-\lambda_i$ are roots of the resultant $R(T)$. Let $\rho(T)=\psi(T^{1/2})$. $\deg_T(\rho)=\displaystyle{m\choose 2}$.}
\end{itemize} 

Therefore the following result holds.
\begin{proposition}\label{prop-roots}
The spread of $\A$ is the biggest real root of $R(T)$. 
The spread of $\A$ is the square root of the biggest real root of $\rho(T)$. The spread of $\A$ is the square root of the biggest real root of the squarefree part of $\rho(T)$. 
\end{proposition}

\section{Reliable computation as proof}\label{sec-reliable}
When a computer does exact integer arithmetic, then so long as \textsl{integer overflow} does not occur, we regard the results as being reliable: as good, in fact, as a mathematical proof.  Given modern computer algebra systems which have what is called ``arbitrary precision'' integer arithmetic, which will only overflow if the entire block of computer memory allocated for the program is filled, we regard this as a solved problem.  Even Python's ordinary integer arithmetic is of this kind.

In this paper we use arbitrary precision integer arithmetic in Maple, which is reliable and reasonably fast.  Maple also implements exact rational arithmetic using arbitrary precision integers.  Our computations usually go nowhere near the memory limits, and we could in fact use short (32 bit) integers most of the time, because the integer coefficients of the characteristic polynomials that we generate are from matrices of small dimension. But because the arbitrary precision arithmetic is fast enough for our purposes, we always use it for simplicity and clarity of programming. We also use interval arithmetic with exact rational bounds when considering roots of characteristic polynomials.

Floating-point arithmetic for the problem we consider here is \textsl{also} reliable in modern computing environments, because rounding errors are quite well understood in these contexts.  In 2023, an excellent interval arithmetic package based on Fredrik Johansson's Arb, with floating point bounds for the intervals, was implemented in Maple~\cite{johansson2017arb}.  Arb has now been merged into \href{https://flintlib.org/doc/}{FLINT}.  Nonetheless we will not use floating-point arithmetic (interval or otherwise) but rather the interval arithmetic with exact rational bounds, which is an older technology (based on Descartes' Rule of Signs~\cite{CollinsAkritas1976,CollinsLoos1976} and implemented in Maple since the 1980s) because it is more like human computation of bounds.

It may seem strange to the readers of this paper that we eschew floating-point arithmetic, even interval arithmetic using floats, in favour of arbitrary-precision rational arithmetic and interval arithmetic using such numbers.  The main drawback of exact rational arithmetic, however, which is that the numbers may grow to be surprisingly large, does not in fact bother us in this application, and the obviation of having to deal with rounding errors is worth the miniscule extra computational cost (for the computations of this paper).  As we will see, the computational cost is dominated by the generation of the matrices to be worked with; the matrices are subsequently processed rapidly using exact integer and symbolic arithmetic to generate characteristic polynomials, resultants, and exact rational intervals containing real roots.  
\section{Brute exhaustive computation for dimension $m=4$}\label{sec-brute-force}
In this section we describe a naive Maple program to exhaustively search for the matrices with maximal spread.  In some sense the computation is not interesting, merely brutal, but in its brutality the results are very clear, as are the potential difficulties in extending this naive method to higher dimension.  
\subsection{The case $a=0$}
There are $2^{m(m+1)/2}$ different $m \times m$ symmetric matrices with two elements.  Therefore there are $2^{10} = 1024$ different four-by-four symmetric matrices with entries either zero or one. We use the \lstinline{CartesianProduct} function from the \lstinline{Iterator} package to generate all of these matrices, and then compute the characteristic polynomials $p(\lambda)$ of each of them.  We then compute the resultant $R(T) = \mathrm{resultant}_\lambda(p(\lambda),p(\lambda+T))$ whose roots are the differences $\lambda_i-\lambda_j$ of the roots of $p(\lambda)$, i.e. of the eigenvalues of $\A$.

There are only $53$ unique resultants generated by this process, in part because many matrices in this collection have the same characteristic polynomial.  One of them, corresponding to the all zero matrix which has spread zero, is discarded.  Now it is a matter of examining each of these resultants to find which has the largest real root (see Prop. \ref{prop-roots}).
 
For instance, one such resultant gives (after computing the square-free factorization to remove multiple roots, and the remaining factor $T$ corresponding to the occurrences of $\lambda_i-\lambda_i$ as root differences)
\begin{equation}
    T^{12}-36 T^{10}+450 T^{8}-2358 T^{6}+5265 T^{4}-4374 T^{2}+729\>.
\end{equation}
We then use the \lstinline{realroot} function on this polynomial, which returns twelve rational intervals containing the root differences $\lambda_i-\lambda_j$.  The intervals do not overlap (by the design of the code) and are sorted in increasing order.  The routine reports that the  largest root $T_s$ of this polynomial, which is the spread of the original matrix, is contained in the interval ${\tfrac{993}{256}} < T_s <  {\tfrac{497}{128}}$, or in decimals using a common shorthand notation, $3.8_{\,\textcolor{red}{789}}^{\,\textcolor{blue}{829}}$ where we have rounded the lower limit down and the upper limit up\footnote{This notation for a narrow interval is (apart from the color choice here, which varies) reasonably widely used, and quite intelligible.  One sees at once an approximation for the true value, and guaranteed bounds for the true value.}.

One matrix representative with this resultant is
\begin{equation}
    \left(\begin{array}{cccc}
1 & 1 & 1 & 1 
\\
 1 & 1 & 1 & 0 
\\
 1 & 1 & 0 & 0 
\\
 1 & 0 & 0 & 0 
\end{array}\right)
\end{equation}
The smallest eigenvalue of this matrix is $-1$, and the largest is approximately $2.87938524157182$ (we can write down an exact expression for the root using in this case the cubic formula because one integer root of the quartic can be exactly factored out, but the complicated resulting expression is not useful).  The point of this example is to show explicitly that the largest root of the resultant does indeed give the spread of the matrix.

Looking over all 52 resultants (ignoring the all-zero matrix) we find that 
\begin{equation}
\left(\begin{array}{cccc}
1 & 1 & 1 & 1 
\\
 1 & 1 & 1 & 1 
\\
 1 & 1 & 1 & 1 
\\
 1 & 1 & 1 & 0 
\end{array}\right)
\end{equation}
gives the largest upper bound on the spread, namely $150163/32768$.  By brute computation, we also see that no other matrix has its upper bound for its largest root larger than the \textsl{lower} bound $300323/65536$ for the largest root of this matrix.

This proves that the matrix above attains the maximal spread in this class of matrices; and this rank-two matrix with a one-by-one block of zeros is exactly of the conjectured form.  This brute computation proves the conjecture is true for $a=0$ and $m=4$.  

The interesting feature is the severe reduction in number of cases that had to be examined, by the fact that only 53 resultants occurred for 1024 different matrices.  We will take systematic advantage of such a reduction when we consider matrices of higher dimension.

But now we do the four-by-four case explicitly with symbolic $a$.
\subsection{Symbolic $a$}
We repeat the computation, exhaustively computing the resultants for all $1024$ distinct matrices.  This time there are only $77$ unique non-trivial resultants.  This is more than the $52$ that occur when $a=0$, but significant reduction has still taken place. 

One instance of a matrix that occurs among these 77 cases is
\begin{equation}
\left(\begin{array}{cccc}
a  & 1 & 1 & 1 
\\
 1 & a  & 1 & 1 
\\
 1 & 1 & 1 & a  
\\
 1 & 1 & a  & a  
\end{array}\right)
    \label{eq:matrixinstance}
\end{equation}
and the corresponding resultant gives, after removing the factor $T^4$,
\begin{align}
& T^{12}+\left(-11 a^{2}+6 a -43\right) T^{10}+\left(41 a^{4}-44 a^{3}+356 a^{2}-232 a +647\right) T^{8}\nonumber\\
&{}+\left(-61 a^{6}+86 a^{5}-1075 a^{4}+1176 a^{3}-3295 a^{2}+3202 a -4129\right) T^{6}\nonumber\\
&{}+\left(31 a^{6}+30 a^{5}+1341 a^{4}+456 a^{3}+6542 a^{2}+5128 a +11560\right) \left(a-1 \right)^{2} T^{4}\nonumber\\
&{}-\left(a^{2}+2 a +46\right) \left(5 a^{4}+8 a^{3}+179 a^{2}+330 a +262\right) \left(a-1 \right)^{4} T^{2}\nonumber\\
&{}+4 \left(5 a^{4}-2 a^{3}+189 a^{2}+404 a +316\right) \left(a-1 \right)^{6}\>.
\label{eq:instance}
\end{align}
Now we must decide if, at any point in the interval $-1 < a < 1$, this polynomial has a root larger than the conjectured maximal spread, which if $-1 < a < -1/3$ is the spread of
\begin{equation}
    \left(\begin{array}{cccc}
a  & a  & 1 & 1 
\\
 a  & a  & 1 & 1 
\\
 1 & 1 & 1 & 1 
\\
 1 & 1 & 1 & 1 
\end{array}\right) \label{eq:twobytwom4}
\end{equation}
while if $-1/3 < a < 1$ the conjectured maximal spread is that of
\begin{equation}
    \left(\begin{array}{cccc}
a  & 1 & 1 & 1 
\\
 1 & 1 & 1 & 1 
\\
 1 & 1 & 1 & 1 
\\
 1 & 1 & 1 & 1 
\end{array}\right)\>.\label{eq:onebyonem4}
\end{equation}
The eigenvalues of the matrix in equation~\eqref{eq:onebyonem4} are
$0$, $0$, $(a+3)/2 \pm \sqrt{a^2-6a+21}/2$ so its spread is $\sqrt{a^2-6a+21}$, while the eigenvalues of the matrix in equation~\eqref{eq:twobytwom4} are $0$, $0$, $1+a \pm \sqrt{a^2-2a+5}$ so its spread is $2\sqrt{a^2-2a+5}$.  We may verify that these spreads are equal if $a=-1/3$.  

\begin{figure}
    \centering
    \includegraphics[width=0.5\linewidth]{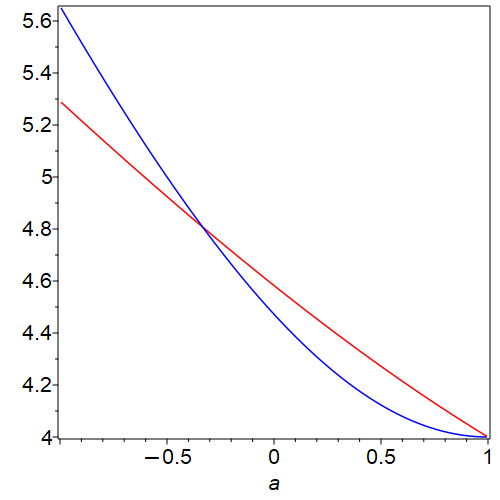}
    \caption{The spread of the rank-two matrix with a two-by-two block of $a$s from equation~\eqref{eq:twobytwom4} is plotted in blue.  The spread of the rank-two matrix with a one-by-one block of $a$s from equation~\eqref{eq:onebyonem4} is plotted in red.  We see that the conjectured maximal spread has a switching point at $a=-1/3$.}
    \label{fig:conjecture4}
\end{figure}

\begin{figure}
    \centering
    \includegraphics[width=0.5\linewidth]{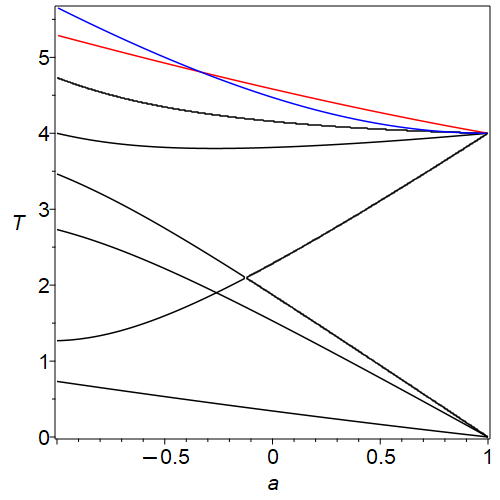}
    \caption{All six positive zeros of the resultant from equation~\eqref{eq:instance} (as functions of the parameter $a$), plotted together with the conjectured maximal spreads from Figure~\ref{fig:conjecture4}. This figure demonstrates (not very convincingly, because the uppermost black curve is quite jagged) that the spread of the matrix in equation~\eqref{eq:matrixinstance} is, for every $a$ in the interval $-1 < a < 1$, smaller than the conjectured maximum (plotted here as the larger of the red and the blue curves).  Notice that at $a=1$ the curves coalesce; but in that case, all the matrices become the unit matrix, which has spread $m$.}
    \label{fig:conjecture4instance}
\end{figure}

We can plot the roots of equation~\eqref{eq:instance} by using \lstinline{implicitplot} from the plots package, and we do so in figure~\ref{fig:conjecture4instance} together with the conjectured bound.  It appears as though all positive roots lie uniformly below the conjectured bound, but it is a bit difficult to see owing to the jaggedness of the uppermost line.  This is a plotting artifact: if we compute the implicit plot on a finer grid, and zoom in, it looks better, though still somewhat hard to interpret near $a=1$ because the largest root of the instance resultant is in fact \textsl{above} the blue curve (though not above the red curve).

To confirm this we may choose a value of $a$, say $a=0.96$ (although to ensure exact arithmetic we write it as $96/100$ or $48/50$ or $24/25$ and Maple will handle the reduction to lowest terms), and compute bounds on the roots of equation~\eqref{eq:instance} by using the \lstinline{realroot} command, this time with a tighter tolerance, $2^{-20}$.  The largest root is contained within the interval
\begin{equation}
    \left[{\frac{137537358207}{34359738368}}, {\frac{17192169777}{
4294967296}}\right]
\end{equation}
which we may write more intelligibly as $4.002863954_{\,\textcolor{red}{72}}^{\textcolor{blue}{\,99}}$.

The value of the blue curve when $a=96/100$ is 
\begin{equation}
    2\sqrt{\frac{2501}{625}} = 4.000799920016\ldots
\end{equation}
which is lower than the lower bound of the instance; but the value of the red curve when $a=96/100$ is
\begin{equation}
    \sqrt{\frac{10101}{625}} = 4.02014\ldots
\end{equation}
which is larger than the upper bound of the instance.

This example shows two things that will be important for our future computations: first, that it is possible for a matrix with rank greater than two to have a spread that is larger than that of a rank-two matrix (indeed, larger than the spread of a matrix which, for other values of $a$, actually attains the maximum spread).  Second, we learn that visual confirmation of maximality is difficult and unconvincing, when the numerical values of the samples are close.

We therefore need some more effective tools than our own eyes to use in deciding with certainty whether or not the root curves ever cross the curves defining the conjectured maximum.  Somewhat whimsically, we called our program to do this our ``Artificial Eye'' or A-eye, for short.  We will develop that in section~\ref{sec:symbolica}.  

\begin{figure}
    \centering
    \includegraphics[width=0.5\linewidth]{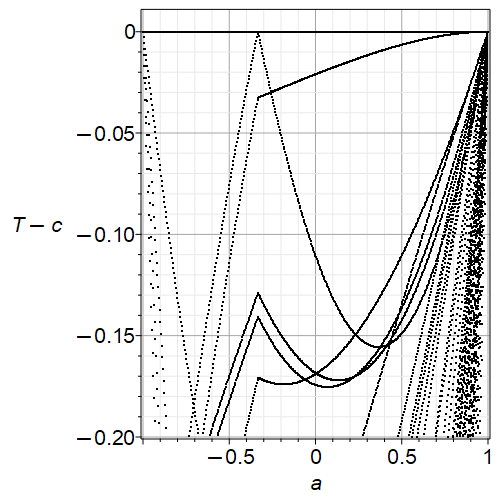}
    \caption{The difference between the computed spread of each of the 77 classes of matrices of dimension $m=4$ and the conjectured maximal spread, at $501$ equally-spaced discrete values of $a$ in the interval $(-1,1)$.    No counterexamples to the conjecture are detected, although the graph is ambiguous near $a=1$ where one of the curves appears to approach the $x$ axis.}
    \label{fig:visualm4}
\end{figure}

For completeness, we computed the maximum zeros of all 77 different nontrivial resultants and subtracted the conjectured maximum (which makes it easier to see differences) and plotted the results. In detail, we sampled each of the 77 different curves at $501$ different values of $a$ across the interval, numerically found the roots in $T$ using \lstinline{fsolve} for each of those values of $a$, and then picked the maximum of the twelve roots.  We subtracted the conjectured maximum.  The result is messy to look at, as you can see in figure~\ref{fig:visualm4}, and as discussed above not completely convincing, but no curves were seen to ever be above zero, even when we zoomed in.  

Again for completeness we give the matrix whose spread came closest to the conjectured maximum:
\begin{equation}
    \left(\begin{array}{cccc}
1 & 1 & 1 & 1 
\\
 1 & a  & 1 & 1 
\\
 1 & 1 & 1 & 1 
\\
 1 & 1 & 1 & a  
\end{array}\right)\>.
\end{equation}
The eigenvalues of this rank 3 matrix are $0$, $a-1$, and $(a+3)/2 \pm \sqrt{a^2-2a+17}/2$. Whatever the value of $a$ in the interval $(-1,1)$, the eigenvalue $(a+3)/2 + \sqrt{a^2-2a+17}/2$ must be the maximum.  The minimum might, a priori, be either $a-1<0$ or $(a+3)/2 - \sqrt{a^2-2a+17}/2$.  By any of several methods we can show that the minimum is actually $a-1$.  Therefore the spread of this matrix is $(1-a)/2 + \sqrt{a^2-2a+17}$.  Comparing this to the conjectured maximal spread we see that it is smaller, intersecting only at $a=1$.  However, it gets very close indeed: the slope of this curve at $a=1$ is exactly the same as the slope of the conjectured maximal spread curve, at $a=1$: both are $-1/2$.  We leave verification of this last fact to the reader.  We have therefore demonstrated by example that matrices of rank other than $2$ can have spread arbitrarily close to the maximal spread.

We do not give the code here that we used to produce this data, because anyone ``skilled in the art'' could reproduce that numerical computation independently.  The graph one produces this way is not a proof (however convincing it is and however finely one samples). We will supply a computational proof in section~\ref{sec:symbolica} using a different method.
\section{Using graph isomorphism to reduce redundant computation}\label{sec-graphs}
Suppose that $\S$ is a member of an unstructured Bohemian family $\Fam$ of matrices.  Then any permutation similarity of $\S$, namely $PSP^{-1}$ where $P$ is a permutation of the identity matrix, is also Bohemian and in the family $\Fam$.  This is also true for some structured Bohemian families, such as symmetric Bohemians, because permutation similarity preserves matrix symmetry.

Now consider the set of all characteristic polynomials $\Char$ of $\Fam$.  Each polynomial $p(\lambda)$ in $\Char$ will have an associated set of matrices in $\Fam$ for which $p(\lambda)$ is the characteristic polynomial.  We may call that set the set of Bohemian companions for $p(\lambda)$.  It is clear that any permutation similarity of a companion for $p(\lambda)$ is still a companion for $p(\lambda)$, so the set of Bohemian companions for $p(\lambda)$ includes all matrices permutation similar to the original choice.  It is, however, possible in general that there may be more than one orbit of permutation similarities in the set of companions for $p(\lambda)$ because having the same characteristic polynomial does not mean matrices are similar, in general.  It is only having the same Frobenius form that means the matrices are similar. But since our family \Fam\ is symmetric, all matrices are diagonalizable and hence if matrices have the same eigenvalues they have the same Frobenius form and thus are similar; so in our case we do not need to worry about this.

It is also possible, however, for two matrices $\A$ and $\B$ in $\Fam$ to be similar but not by a permutation similarity.  This is unlikely, but might happen.  What this means in practice is that the characteristic polynomial would be discovered at least twice even in an efficient exhaustive search, such as we are about to describe. This double-counting possibility will not, therefore, cause us to miss any candidates for maximal spread, but only mean that the search was somewhat less efficient than it might have been.

Now we consider graph isomorphism.  If we can express $\S$ in some form that is preserved under permutation similarity but contains a matrix that can be interpreted as an adjacency matrix for a graph, then we will be able to identify a canonical member of the orbit of $\S$ under permutation similarity.  Working only with that canonical member will reduce the amount of computation that we have to do.

In the case of symmetric matrices whose elements are either $0$ or $1$, there is a natural identification with adjacency matrices; except that because the diagonal elements may be $1$, we have to allow the possibility of self-loops in the graphs.

In the case of symmetric matrices whose elements are either $1$ or $a$ where $a$ is not necessarily $0$, we may write
\begin{equation}\label{eq:graphdecom}
  \S = \mathbf{E} + (a-1)\G
\end{equation}
where $\mathbf{E} = \mathbf{e}\mathbf{e}^T$ is the matrix of all $1$s, $\mathbf{e}$ is the column vector of all $1$s, and $\G$ is a symmetric $0$--$1$ matrix, which we may interpret as an adjacency matrix for a graph, again with the possibility of self-loops.

Every graph isomorphism can be represented as a permutation similarity of~$\G$.

If we can identify a canonical representation of the graph (and therefore a canonical adjacency matrix $\G$) then we will have found a canonical representation for $\S$.

If we find \textsl{all} canonical matrices $\S$ in the family $\Fam$, then we will have found all possible characteristic polynomials.  It is possible that there will still be some repeat polynomials in this list, because as mentioned previously some matrices not similar by permutation  may still have the same characteristic polynomial, but there will be many fewer than naively enumerating all matrices.

This suggests that we enumerate canonical looped graphs and their matrices $\G$ using the \texttt{nauty} package~\cite{MCKAY201494} (which is partially implemented in the \texttt{GraphTheory} package in Maple) and use those to build matrices $\S$ which we take the characteristic polynomials of (or compute eigenvalues of, if numerical values of $a$ are known).

Comparing the number of symmetric matrices with population $\{1,a\}$ at each dimension and the number $U_m$ of nonisomorphic unlooped graphs at each dimension $m$ (this is A000088 in the OEIS) together with the number $L_m$ of nonisomorphic graphs with loops, which is A000666 in the OEIS, as in Table~\ref{tab:nonisomorphic}, we see a very substantial reduction in the number of matrices that need to be processed.

\begin{table}
  \centering
  \caption{Reduction by canonicalization. The first column is the dimension, $m$. The second column is the number of different $m\times m$ symmetric matrices with each entry being either $1$ or $a$.  The third column is the number $U_m$ of canonical graphs.  The third column is the number of semi-canonical matrices with potentially nonzero diagonals allowing self loops. 
 The fourth column contains entries from the OEIS sequence A000666, namely the numbers $L_m$ of nonisomorphic graphs with loops. }\label{tab:nonisomorphic}
  \begin{tabular}{c|c|c|c|c}
    \toprule
    $m$ & $2^{m(m+1)/2}$& $U_m$ & $(2^m-1)U_m$ & $L_m$ \\
    \midrule
 2 & 8 & 2 & 6 & 6 
\\
 3 & 64 & 4 & 28 & 20 
\\
 4 & 1.024 & 11 & 165 & 90 
\\
 5 & 32.768 & 34 & 1.054 & 544 
\\
 6 & 2.097.152 & 156 & 9.828 & 5.096 
\\
 7 & 268.435.456 & 1.044 & 132.588& 79.624 
\\
 8 & 68.719.476.736 & 12.346 & 3.148.230& 2.208.612 \\
 9 & 35.184.372.088.832 & 274.668 & 140.355.348& 113.743.760 \\
\bottomrule
\end{tabular}

\end{table}
The fifth column is because the non-isomorphic graphs generated by \lstinline{geng} are all free of self-loops, so the diagonal entries must be zero.  But in
our application the diagonal entries of $\S$ are either $1$ or $a$.  To allow diagonal entries to be $1$, we must allow self-loops.  One way to do this is to colour the vertices, say blue or yellow.  A blue vertex will correspond to a self loop being present, while a yellow vertex has no self loop.  Then if we postprocess the output of \lstinline{geng} using \lstinline{vcolg} we can generate all nonisomorphic graphs with self loops.

In practice what we have actually done (so far) is instead to use the Maple version, \lstinline{GraphTheory:-NonIsomorphicGraphs} and then account for all possible diagonals ourselves, discarding one of a pair that gives the same spread; this generates $(2^m-1)U_m$ different graphs.  This is reported in the next-to-last column.

The spread of the eigenvalues if the diagonal entries are all $1$ is the same as the spread of the eigenvalues if the diagonal entries are all $a$, so we only need one of those two; but in
order to account for all the possibilities (with the simple implementation of \lstinline{geng} found in Maple) we need to add all $2^m-1$ possible choices of diagonals. This still reduces the effort overall, as seen in the third column of counts.

\subsection{Formalizing the argument}
Fix the dimension $m$.  Let $\mathbb{S}$ represent the family of all symmetric matrices of dimension $m$ whose entries are either $1$ or the symbol $a$.  We interpret $a$ to be a fixed real number in the semiclosed interval $-1 \le a < 1$; \textcolor{black}{that is, $\mathbb{S}=S_m(\{1,a\})$}.  Let $\mathbb{G}_U$ represent the family of adjacency matrices of undirected graphs on $m$ vertices; each entry of an adjacency matrix $\G \in \mathbb{G}_U$ is either zero or one, and moreover the diagonal entries of $\G$ are all zero.


Then let $\mathbb{D}$ be the set of all diagonal matrices (of dimension $m$) whose diagonal entries are either $1$ or $0$.  Later, we will remove the diagonal whose entries are all $1$, and keep only the one whose entries are all $0$, but for simplicity now we include both.  The cardinality of $\mathbb{D}$ is thus $2^m$.

Let $E = \mathbf{e} \mathbf{e}^T$ be the rank-one matrix of all ones, formed by the outer product of the column vector $\mathbf{e}$ of all ones with itself. 

Then we have the following propositions, given without proof. 
\begin{enumerate}
\item Every symmetric matrix $\S \in \mathbb{S}$ can be uniquely represented as
\begin{equation}
    \S = \E + (a-1)(\D + \G) 
\end{equation}
using a unique matrix $\D \in \mathbb{D}$ and a unique adjacency matrix $\G \in \mathbb{G}_U$. Indeed the correspondence is one-to-one between $\mathbb{S}$ and the direct sum $\mathbb{D} \oplus \mathbb{G}_U$.
\item The orbits $\P \S \P^T$ of $\S \in \mathbb{S}$ under similarity transforms by \textsl{permutation matrices} $P$ (for which $P^{-1} = P^T$) are disjoint in $\mathbb{S}$ and clearly cover $\mathbb{S}$.
\item $\P\S\P^T = \E + (a-1)(\P\D\P^T + \P\G\P^T)$ so $\S$ is permutation similar to (at least) one pair in $\mathbb{D} \oplus \mathbb{G}_U$.
\item The graph with adjacency matrix $\P \G \P^T$ is isomorphic to the graph with adjacency matrix $\G$.
\item Given a canonical representative adjacency matrix $\G$ for all graphs isomorphic to the graph with adjacency matrix $\G$, computing characteristic polynomials of $\E + (a-1)(\D_k + \G)$ for all $D_k \in \mathbb{D}$ will result in all possible characteristic polynomials in the orbit of $\S$.
\item The eigenvalues of $\S_1 = \E + (a-1)(\D_1 + \G)$ and $\S_2 = \E + (a-1)((\D_1 + \I) + \G)$ are related; if the eigenvalues of $\S_1$ are $\lambda_j$ for $1 \le j \le m$ then the eigenvalues of $\S_2$ are $\lambda_j + 1$ for $1 \le j \le m$.  This means that the \textsl{spread} of the two matrices, i.e. the difference between the largest eigenvalue $\lambda_1$ and the smallest eigenvalue $\lambda_m$ (since the eigenvalues are real, we may assume that they are sorted so $\lambda_m \le \lambda_{m-1} \le \cdots \le \lambda_2 \le \lambda_1$) is the same for both matrices.  
\end{enumerate}

\section{The special case $a=0$}\label{sec-casezero}
\textcolor{black}{In this section, we prove that the conjecture is true for $2\leq m\leq 8$ and $a=0$.}

\begin{theorem}
If $a=0$ and $m \in \{2, 3, 4, 5, 6, 7, 8\}$, then the maximal spread of $\S \in S_m(\{1,a\})$ is attained by a rank-two matrix $\S^* \in S_m(\{0,1\})$, and is explicitly given by 
    $\sqrt{(m+k)^2- 4k^2}$ where $k = [ m/3 ]$.  That is, $k$ is the nearest integer to $m/3$.
    The matrix $\S^*$ is, as conjectured by Fallat and Xing, any of those permutationally equivalent to the one of the form
    \textcolor{black}{\begin{equation}
\left(\begin{array}{cc}
0\,\J_k&\J_{k,m-k}\\
\J_{m-k,k}&\J_{m-k,m-k}
\end{array}\right)
    \end{equation}}
\end{theorem}
The proof is computational. The algorithm is very simple, and can be described in words.  We choose a dimension $m \ge 3$ (the case $m=2$ is simple enough to be carried out using exhaustive computation even by hand).  We set $a=0$.  We use the nauty package~\cite{MCKAY201494} as implemented in the Maple \lstinline{GraphTheory} package to generate all the nonisomorphic graphs on $m$ vertices.  For each of those graphs, we generate the corresponding possible matrices by walking through every possible diagonal (except one: the matrix with diagonal $[1,1,\ldots,1]$ will have the same spread as the matrix with diagonal $[a,a,\ldots,a]$ and so we only look at one of the two).

We compute the rank of this constructed matrix.  If the rank is two, we discard the matrix, because Fallat and Xing have already proved that the matrices of the form above indeed have the stated maximal spread, over all rank-two symmetric matrices with entries either $1$ or $a$.

We then compute the Mirsky bound for this particular matrix (see \eqref{MirskyBound}).  If this bound is strictly smaller than the spread of the conjectured optimal rank two matrix, which is $\sqrt{(m+k)^2-4k^2}$, then we discard the matrix and proceed to the next one.  This step potentially saves some computational effort.  

If the Mirsky bound of our generated matrix is greater than or equal to the spread of the conjectured optimal rank two matrix, then the matrix is a potential counterexample to the conjecture.  We then proceed with the following more expensive steps.

We compute the characteristic polynomial $p(\lambda)$ using exact integer arithmetic.  We then compute a shifted polynomial $q(\lambda,T) = p(\lambda+T)$ in a new symbolic variable $T$.  The resultant $R(T) = \mathrm{resultant}_{\lambda}(p(\lambda),q(\lambda,T)$ then has the following properties (see Subsection \ref{subsection-resultant}).
\begin{enumerate}
    \item The roots of $R(T)$ are of the form $\lambda_i-\lambda_j$ where $\lambda_i$ and $\lambda_j$ range over all possible roots of $p(\lambda)$.
    \item $T^m | R(T)$ because each $\lambda_i - \lambda_i$ is a root of $R(T)$
    \item $\rho = R(T)T^{-m}$ is a function of $T^2$ because both $\lambda_i-\lambda_j$ and $\lambda_j-\lambda_i$ are roots of $R(T)$.
    \item The largest root of $\rho(\tau)$, where $\tau=T^2$, is the square of the spread of the matrix with characteristic polynomial $p(\lambda)$.
\end{enumerate}
We therefore construct the squarefree part of $\rho(\tau)$ where $\tau = T^2$ and find its largest root.  We use \lstinline{realroot} to do so, which is based on Descartes rule of signs~\cite{CollinsAkritas1976}.  In fact, \lstinline{realroot} computes the square free part of its input automatically, so we do not have to explicitly carry out that step.  The routine \lstinline{realroot} works in exact rational interval arithmetic and returns  intervals guaranteed to contain each real root.  Although it is an old program, and the method it uses was published fifty years ago, it remains competitively efficient, and is appropriate for this application.

Our overall program returns every polynomial (and a matrix corresponding to that polynomial) with an upper bound of a root larger than a factor $999/1000$ of the conjectured maximal spread.  This program is therefore guaranteed to detect any matrices with spread equal to or larger than the conjectured maximal spread, possibly together with some false positives where the bound was larger but the spread was not. Each returned polynomial beyond the resultant of the rank two matrix would need to be separately examined.

In practice this was simple.  Our program reported that for dimensions $m=3$, $4$, $5$, $6$, $7$, and $8$, 
no counterexamples were generated in any of our runs. Therefore, no matrices exist at these dimensions with larger spread than the rank two matrix which was conjectured to have maximal spread.  This exhaustively proves that the conjecture is true for these dimensions.


The computing times for the various dimensions were, on a 2017 Microsoft Surface Pro using Windows 10 and running Maple 2024, as reported in table~\ref{tab:computingtimesazero}.
\begin{table}[t]
    \centering
    \begin{tabularx}{0.8\textwidth}{ >{\centering\arraybackslash}X| >{\centering\arraybackslash}X}
    \toprule
     $m$ & time (s) \\
     \midrule
       $2$  & $0.11$ \\
       $3$  & $0.25$ \\
       $4$  & $0.392$ \\
       $5$  & $1.375$ \\
       $6$ & $9.3$\\
       $7$ & $137$ ($2$m $17$s)\\
       $8$ & $4300$ ($1.2$ hr) \\
       \bottomrule
    \end{tabularx}
    \caption{Computing times to search all canonical symmetric $\{0,1\}$ matrices to verify that the conjectured maximal spread was correct. Growth is (as expected) faster than exponential in $m$.  The case $m=9$ is estimated to take approximately a week on this machine, a 2017 vintage Microsoft Surface Pro.}
    \label{tab:computingtimesazero}
\end{table}
\section{Symbolic $a$}\label{sec:symbolica}
\textcolor{black}{In this section it is shown that the conjecture is valid for $m \in \{2, 3, 4, 5, 6, 7\}$ and $a$ within the interval $(-1,1).$}

\begin{theorem}
If $-1 < a < 1$ and $m \in \{2, 3, 4, 5, 6, 7\}$, then the maximal spread of $\S \in S_m(\{a,1\})$ is attained by a rank-two matrix $\S^* \in S_m(\{a,1\})$, and is explicitly given below, where $k = [ m/(3+a) ]$.  
    The matrix $S^*$ is, as conjectured by Fallat and Xing, any of those permutationally equivalent to the one of the form
    \begin{equation}
     \textcolor{black}{  \left(\begin{array}{cc}
a\,\J_k&\J_{k,m-k}\\
\J_{m-k,k}&\J_{m-k,m-k}
\end{array}\right)}
    \end{equation}
Moreover, the interval $-1 < a < 1$ is divided by the breakpoints $\alpha_{m,\ell} = 2m/(2\ell +1)-3$ for $\ell = \ceil{(m-2)/4}$ up to $\floor{(m-1)/2}$.  Over intervals defined by these breakpoints, the optimal block size $k$ is constant. For instance, if $m=4$ then the breakpoints predicted by that formula are $-1$ (which we discard as redundant) and $-1/3$, so the two intervals are $-1 < a < -1/3$ on which, by $k=[m/(a+3)]$, we have $k=2$ and the other interval $-1/3 < a < 1$ on which $k=1$.  At the breakpoint $a=-1/3$, both block sizes give the same maximal spread.  

The formula given above in terms of $\ell$ and the floor functions will give a redundant $-1$ when $m$ is odd and a redundant $1$ when $m = 2 \mod 4$.
\end{theorem}

Again the proof is computational, but this time because we carry a symbol $a$ throughout the computation the cost in memory and time is larger, and on this small machine we were able only to go to dimension $m=7$.

The computing time for $m=3$ was $0.3$ seconds.  The computing time for $m=4$ was 2 seconds.  The computing time for $m=5$ was $31$ seconds. The computing time for $m=6$ was $16.25$ minutes. The computing time for dimension $m=7$ was approximately 36 hours.  We estimate that $m=8$ would take over forty days, if it could be completed without memory issues on this small machine.

For each dimension, the script examined all possible semi-canonical matrices in turn.  The script computed the characteristic polynomial $p(\lambda; a)$ for each matrix, and thence the resultant $R(T;a) = \mathrm{resultant}_\lambda(p(\lambda;a),p(\lambda+T;a))$. Since $R(T;a)$ has a zero of multiplicity $m$ at $T=0$ the script divides this out.  The remaining polynomial is a function of $\tau = T^2$ because not only is $T=\lambda_i-\lambda_j$ a zero of $R(T;a)$, so is $-T=\lambda_j-\lambda_i$. The script factors this bivariate polynomial, over the integers, to make subsequent resultant computations faster.  The rank-two maximum was written as the polynomial $U(\tau;a) = \tau - (a^2 + 2a - 3)k^2 + 2m(1 - a)k + m^2$, where $k=[m/(a+3)]$ is a constant integer on each appropriate subinterval of $-1 < a < 1$.  
The script computes $Z(a) = \mathrm{resultant}_\tau( \rho(\tau;a), U(\tau;a))$, which is a polynomial in $a$.  If this is identically zero, then we have found a matrix equivalent to the rank-two matrix in that it has the same spread.  Finally, the script computes the Sturm sequence of $Z(a)$ on the interval where $k$ is constant.  If this is nonzero then $\rho(\tau;a)$ intersects $U(\tau;a)$ in this interval.  This will always happen on the final interval where $a=1$ so that case is discarded.  The script found no nontrivial intersections.  

Since a separate computation established that for $a=0$ the values $\rho(\tau;0)$ were always smaller than $U(\tau;0)$, by continuity we conclude that no counterexamples were found. Indeed this script eliminated all possible counterexamples and proved the theorem for $m=3$, $4$, $5$, $6$, and $7$.

\noindent \textbf{Remark}. This ``exhaustive computational proof'' relies on several subroutines supplied by Maple.  These include \lstinline{factor}, \lstinline{resultant}, \lstinline{sturm}, the routines from the LinearAlgebra package, and those from the GraphTheory package. With the exception of the GraphTheory package, which was revised within the last ten years, all of these subroutines and packages have been tested extremely intensively by a great many people and in a great number of applications.  Our level of confidence in them is, in fact, greater than our level of confidence in the typical long ``by hand'' proof.  We acknowledge that not everyone will agree.
As for the GraphTheory package, it is really an interface to the freely-available \lstinline{nauty} package~\cite{MCKAY201494}, which again has been very widely used and tested.  Again, we have very high confidence in its correctness. \textcolor{black}{\textbf{ Our own code, calling all of these packages, is relatively minimal and is available for inspection at \url{https://github.com/rcorless/BohemianSpread} on {Rob Corless' GitHub page}} }.

In this section we show how to construct, in effect, a symmetric companion matrix for the resultant polynomial $\res(T)$, which we will call $\Sc$.  Then the largest root of the resultant, which gives the spread, will be the largest eigenvalue of~\Sc.

It may seem quixotic to start with a matrix eigenproblem, move to a resultant polynomial, and then move back to a different eigenproblem, but there is some benefit to this as we will see.  One difference is that if the original symmetric matrix $\S$ is of dimension $m \times m$, then the symmetric companion matrix will be of dimension $m^2 \times m^2$, which does not seem like progress.

Note first that the resultant has $m^2$ roots $T_{i,j} = \lambda_i - \lambda_j$ where the $\lambda_i$ are the $m$ real eigenvalues of $\S$.  There are $m^2$ of these differences, $m$ of which are zero.  Therefore no matter which companion we choose (there are infinitely many to choose from) the dimension will always be $m^2 \times m^2$.  It might be possible to remove the $m$ zero eigenvalues by (for example) deflation, but we will not pursue that here.

\section{Another Approach}\label{sec:another}

In this section we show how to construct, in effect, a symmetric companion matrix for the resultant polynomial $\res(T)$, which we will call $\Sc$.  Then the largest root of the resultant, which gives the spread, will be the largest eigenvalue of~\Sc.

It may seem quixotic to start with a matrix eigenproblem, move to a resultant polynomial, and then move back to a different eigenproblem, but there is some benefit to this as we will see.  One difference is that if the original symmetric matrix $\S$ is of dimension $m \times m$, then the symmetric companion matrix will be of dimension $m^2 \times m^2$, which does not seem like progress.

Note first that the resultant has $m^2$ roots $T_{i,j} = \lambda_i - \lambda_j$ where the $\lambda_i$ are the $m$ real eigenvalues of $\S$.  There are $m^2$ of these differences, $m$ of which are zero.  Therefore no matter which companion we choose (there are infinitely many to choose from) the dimension will always be $m^2 \times m^2$.  It might be possible to remove the $m$ zero eigenvalues by (for example) deflation, but we will not pursue that here.

One approach to finding a companion is to build it from the symmetric B\'ezout matrix, which gives a matrix polynomial in $T$ of dimension $m$ and degree $m$ with symmetric coefficients.  There are several methods to construct a symmetric linearization of this matrix polynomial, such as described in~\cite{amiraslani2009linearization} or in~\cite{higham2007symmetric}.  However, there is a simpler and more direct method, which uses some properties of tensor products that we list below.  These ``Facts'' are taken from~\cite{reams2013partitioned}, and we assume that all dimensions are compatible.

\begin{enumerate}
  \item Fact 2: $(\A + \B)\otimes \C = \A\otimes \C + \B\otimes \C$
  \item Fact 3: $\A\otimes(\B + \C) = \A\otimes \B + \A\otimes \C$
  \item Fact 5: $(\A\otimes \B)^T = \A^T \otimes \B^T$
  \item Fact 8: $(\A\otimes \B)(\C\otimes \D) = (\A\B)\otimes (\C\D)$
\end{enumerate}
Given these facts, it is straightforward to prove that the matrix
\begin{equation}\label{eq:resultantcompanion}
  \Sc := \mathbf{S}_m \otimes \mathbf{I}_m - \mathbf{I}_m \otimes \mathbf{S}_m
\end{equation}
has eigenvalues $\lambda_i - \lambda_j$ with eigenvectors $\mathbf{u}_i \otimes \mathbf{u}_j$ where the vectors $\mathbf{u}_i$ are themselves the eigenvectors of $S_N$ associated with the eigenvalues $\lambda_i$.

The spread of the matrix $\mathbf{S}_m$ is therefore given by the largest eigenvalue of the symmetric matrix $\Sc$.  Suppose the eigenvalues of $S_m$ are sorted $\lambda_1 \ge \lambda_2 \ge \cdots \ge \lambda_m$.  Then the two largest magnitude eigenvalues of $\Sc$ are $\lambda_1 - \lambda_m$ and its negative, $\lambda_m - \lambda_1$.  These eigenvalues might, a priori, be multiple.

We may use one more fact here, namely the Rayleigh quotient, to give lower bounds on the maximum spread.  Suppose we have two random vectors $\mathbf{x}$ and $\mathbf{y}$.  Then
\begin{equation}\label{eq:Rayleighbound}
  \frac{\mathbf{x}^T\otimes\mathbf{y}^T \Sc \mathbf{x}\otimes\mathbf{y}}{\mathbf{x}^T\mathbf{x}\mathbf{y}^T\mathbf{y}} \le \lambda_1 - \lambda_N
\end{equation}
because the largest eigenvalue of a symmetric matrix maximizes the Rayleigh quotient.  This suggests that we can get lower bounds for the spread of a matrix by making good choices for $\mathbf{x}$ and $\mathbf{y}$.

This can be done efficiently by making $\mathbf{x}$ a good approximation for the eigenvector $\mathbf{u}_1$ corresponding to $\lambda_1$ and $\mathbf{y}$ a good approximation for the eigenvector $\mathbf{u}_N$.

{\color{black} \section{Concluding remarks}\label{sec-conclusions}
By combining exact arithmetic, graph-theoretic reductions, and algorithmic techniques from symbolic computation,
we have provided new computational evidence in support of the \textsl{2012 Fallat and Xing  Conjecture} on the spread of symmetric matrices. Our results confirm the conjecture in several cases, including $S_m([0,1])$ for $m$ divisible by $3$,
and for $2\leq m \leq 8$. Additionally, we prove the conjecture for $S_m([a,1])$
with  $-1 < a < 1$ and $m \in \{2,3,4,5,6,7\}$.  

\textcolor{black}{We did attempt to prove the conjecture in general, and at one point we thought we had succeeded.  The Mirsky bound is attained when the arithmetic mean of the extremal eigenvalues is itself a multiple eigenvalue (with multiplicity $m-2$), and we found a way to use matrix cubic surds to make this happen; however, there was a stubborn ``gap'' in our purported general proof, which we were unable to close.  It is interesting that this ``multiple-root average'' characterization of the maximum spread is somehow ``almost true'' for the Fallat and Xing conjecture.  The rank-two matrix has two extremal eigenvalues, and all other roots are multiple (being zero), but zero is not the average of the extremal eigenvalues.  We still feel, in spite of our experience, that an argument along the lines of what we tried will ultimately succeed.}
}

\section*{Acknowledgments}
\textcolor{black}{We would like to thank Jane Breen for pointing out Theorem 5 of \cite{breen2022maximum}, during the conference \href{https://cs.uwaterloo.ca/~lrafiees/Lalo60.html}{Matrices and Polynomials in Computer Algebra}, allowing us to prove the 2012 Fallat and Xing  Conjecture for $S_{3q}([0,1])$ (see Subsection \ref{subsec-multiple}).} We also remember with gratitude Nick Higham's early and continued interest in Bohemian Matrices. See his column \url{https://www.siam.org/publications/siam-news/articles/rhapsodizing-about-bohemian-matrices/} in SIAM News 2018.
The support of the Rotman Institute of Philosophy is gratefully acknowledged.

\section*{Funding}
This work was partially supported by NSERC grant RGPIN-2020-06438.


\end{document}